\documentclass[12pt]{amsart}
 \usepackage{graphicx}
\usepackage{color}
\usepackage{amsmath,graphics}
\usepackage{amsfonts,amssymb,xypic}

\theoremstyle{plain}
\newtheorem*{theorem*}{Theorem}
\newtheorem*{lemma*} {Lemma}
\newtheorem*{corollary*} {Corollary}
\newtheorem*{proposition*} {Proposition}
\newtheorem{theorem}{Theorem}[section]

\theoremstyle{remark}

\newtheorem*{claim}{Claim}

\theoremstyle{definition}

\def\tpm{[t^{\pm 1}]}
\def\rt{R\tpm}

\def \Z {\Bbb{Z}}

\def\op{\operatorname}

\def\det{\op{det}}

\def\gl{\mbox{GL}}

\def\id{\op{id}}

\def\Z{\Bbb{Z}}

\def\N{\Bbb{N}}

\def\part{\partial}

\def\a{\alpha}
\def\g{\gamma}

\def\bp{\begin{pmatrix}}

\def\sm{\setminus}

\def\ep{\end{pmatrix}}
\def\bn{\begin{enumerate}}

\def\en{\end{enumerate}}
\def\ba{\begin{array}}
\def\ea{\end{array}}

\def\S{\Sigma}

\def\a{\alpha}

\def\ti{\tilde}

\def\fr12{\frac{1}{2}}

\def\be{\begin{equation} }
 \def\ee{\end{equation}}

\def\G{\Gamma}

\def\ol{\overline}

\begin{document}

\title{Twisted Reidemeister torsion, the Thurston norm and fibered manifolds}

\author{Stefan Friedl}
\address{Mathematisches Institut\\ Universit\"at zu K\"oln\\   Germany}
\email{sfriedl@gmail.com}


\date{\today}

\begin{abstract}
We prove  that the twisted Reidemeister torsion of a 3-manifold corresponding to a fibered class is monic
and we show that it gives lower bounds on the Thurston norm.
The former fixes a flawed proof in \cite{FV10}, the latter gives a quick alternative argument for the main theorem of \cite{FK06}.
 \end{abstract}

\maketitle

\section{Introduction}

A \emph{$3$-manifold pair} is a pair $(N,\phi)$ which   consists of an orientable, connected, compact $3$-manifold $N$ with empty or toroidal boundary
and a primitive class  $\phi\in H^1(N;\mathbb{Z}) = \mbox{Hom}(\pi_1(N),\mathbb{Z})$.
 The \emph{Thurston norm} (see \cite{Th86}) of $\phi$ is defined as
 \[
\|\phi\|_{T}=\min \{ \chi_-(\S)\, | \, \S \subset N \mbox{ properly embedded surface dual to }\phi\}.
\]
Here, given a surface $\S$ with connected components $\S_1\cup\dots \cup \S_k$, we define
\[ \chi_-(\S)=\sum_{i=1}^k \max\{-\chi(\S_i),0\}.\]
We furthermore say that $\phi$ is a fibered class if there exists
 a fibration $p\colon N\to S^1$ such that the induced map $p_*\colon \pi_1(N)\to \pi_1(S^1)=\mathbb{Z}$ coincides with $\phi$.

The Thurston norm can be viewed as a generalization of the genus of a knot and fibered classes
are a generalization of fibered knots. It is well known that the Alexander polynomial of a knot contains information about the knot genus and about fiberedness.

This relationship has been generalized lately to twisted invariants.
Recall that given a $3$-manifold pair $(N,\phi)$ and a representation $\a\colon \pi_1(N)\to \gl(k,R)$ over a domain $R$
we can  consider the twisted Reidemeister torsion $\tau(N,\phi\otimes \a)\in Q(t)$,
where $Q$ is the quotient field of $R$.
Note that $\tau(N,\phi\otimes \a)\in Q(t)$ is well-defined up to multiplication by an element in $Q(t)$ of the form
$\pm rt^k$ where $r\in \det(\a(\pi_1(N)))$ and $k\in \Z$.
We refer to Sections \ref{section:deftorsion} and \ref{section:deftorsion2} and to \cite{FV10} for details.
This invariant can be viewed as the generalization of the Alexander polynomial of a knot and it is closely related to the twisted Alexander polynomials of Lin \cite{Li01} and Wada \cite{Wa94}. See \cite{Ki96} and \cite{FV10} for details.

Given
\[ f(t)=a_rt^r+a_{r+1}t^{r+1}+\dots +a_st^s\in \rt\]
with $a_r,a_s\ne 0$ we define $\deg(f(t))=s-r$.
We furthermore say $f(t)$ is \emph{monic} if $a_r$ and $a_s$ are equal to $\pm 1$.
Given $f(t)=p(t)/q(t)\in Q(t)\sm \{0\}$ we define
\[ \deg(f(t))=\deg(p(t))-\deg(q(t)).\]
We say $f(t)\in Q(t)$ is monic if it is the quotient of two monic polynomials in $\rt$.

We start out with the following result.

\begin{theorem} \label{mainthm1}
Let $(N,\phi)$ be a fibered 3--manifold pair with
 $N\ne S^1\times D^2$ and $N\ne S^1\times S^2$. Let
${\a}\colon \pi_1(N)\to \mbox{GL}(k,R)$ be a  representation. Then  $\tau(N,\phi\otimes \a) \in Q(t)$ is monic
and we have
\[ \deg(\tau(N,\phi\otimes \a))=k \cdot \|\phi\|_T.\]
\end{theorem}

The following special cases have been proved before:
\bn
\item Cha \cite{Ch03} showed that twisted Alexander polynomials (which have in general a much larger indeterminacy) of fibered knots are monic.
\item Goda, Kitano and Morifuji \cite{GKM05} showed that the twisted  Reidemeister torsion of a fibered knot is monic,
the same proof also works for any fibered $3$-manifold with non-trivial boundary.
\item In \cite{FK06} it is shown that twisted Alexander polynomials of fibered 3-manifolds are monic.
\en
None of the above  proofs  can be extended in a clear way to provide a proof of Theorem \ref{mainthm1}.
Theorem \ref{mainthm1} was given in \cite{FV10} and a sketch of a short proof was given.
Unfortunately the sketch was too simple-minded and we now give a correct proof of this result.

In this paper we also give a quick proof of the following theorem which was first obtained in \cite{FK06}.

\begin{theorem}\label{mainthm2}
Let $(N,\phi)$ be a 3--manifold pair
and let $\a\colon \pi_1(N)\to \gl(k,R)$ be a representation over a domain $R$. If $\tau(N,\phi\otimes \a)\ne 0$, then
\[ \deg(\tau(N,\phi\otimes \a))\leq k \cdot \|\phi\|_T.\]
\end{theorem}

The proof in \cite{FK06}, as basically all proofs relating (twisted) Alexander polynomials
 to the knot genus and the Thurston norm, relies on a Mayer--Vietoris sequence which relates  the (twisted) Alexander module
 to the homology of a Thurston norm minimizing surface. The proof we give in this paper is quite different.
 It uses an appropriately chosen CW-complex structure for $N$ to calculate the twisted Reidemeister torsion.
 This approach allows us to give a proof which is considerably shorter than the proof in \cite{FK06}.
  Our proof of Theorem \ref{mainthm2} can also be easily adapted to give alternative proofs of \cite[Theorem~1]{Tu02}, \cite[Theorem~10.1]{Ha05}
and \cite[Theorem~1.2]{Fr07}.

Note that a converse to Theorem \ref{mainthm2} was proved in \cite{FV12a} (or alternatively, see \cite[Theorem~1]{FV08} combined with \cite[Theorem~1.1]{PW12}),
namely given a non-fibered class $\phi$ there exists a representation $\a$ such that $\tau(N,\phi\otimes \a)$ is zero.
(See also \cite{FV11} for an earlier and weaker result.)
Also, in \cite{FV12b} it is shown that if $N$ is  irreducible and not a closed graph manifold, e.g. if $N$ is hyperbolic,
then there exists a representation such that the twisted Reidemeister torsion detects the Thurston norm of a given $\phi$.

\subsection*{Convention.} All manifolds are assumed to be connected, compact and oriented, unless it says specifically otherwise.

\subsection*{Acknowledgment.} We are grateful to Wolfgang L\"uck for pointing out the  flawed argument in \cite{FV10}.
We also wish to thank J\'er\^{o}me Dubois for helpful conversations.

\section{Definitions and Preliminaries}

\subsection{Definition of twisted Reidemeister torsion}\label{section:deftorsion}

Let $N$ be a 3-manifold and $X\subset N$ a subspace.
We write $\pi=\pi_1(N)$. Let $\gamma\colon \pi\to \gl(k,Q)$ be a representation over a field $Q$.
We endow $N$ with a finite CW--structure such that $X$ is a subcomplex.
We denote the universal cover of $N$ by $p\colon \tilde{N} \to N$ and we write $\ti{X}:=p^{-1}(X)$.
Recall that there exists a canonical left $\pi$--action on the universal cover $\tilde{N}$ given by deck transformations. We
consider the cellular chain complex $C_*(\tilde{N},\ti{X})$ as a right $\mathbb{Z}[\pi]$-module by defining $\sigma \cdot
g\mathrel{\mathop:}= g^{-1}\sigma$ for a chain $\sigma$ and some $g\in \pi$.

Using the representation $\g$ we can  view $Q^k$ as a left module over $\Z[\pi]$. We can therefore consider
the  $Q$--complex
\[ C_*(\tilde{N},\ti{X})\otimes_{\mathbb{Z}[\pi]}Q^k.\]
We now endow the free $\mathbb{Z}[\pi]$--modules $C_*(\tilde{N},\ti{X})$
with a basis by picking lifts of the cells of $N\sm X$ to $\tilde{N}$.
Together with the canonical basis $v_1,\dots,v_k$ for $Q^k$ we can now view  the $Q$--complex $C_*(\tilde{N},\ti{X})\otimes_{\mathbb{Z}[\pi]}Q^k$
as a complex of based $Q$--vector spaces.

If this complex is not acyclic, then we define $\tau(N,X,\gamma)=0$.
Otherwise we denote by $\tau(N,X,\gamma)\in Q\setminus \{0\}$ the  Reidemeister torsion of this based $Q$--complex.
We will not recall the definition of Reidemeister torsion, referring instead to the many excellent expositions,
 e.g. \cite{Mi66} and \cite{Tu86,Tu01}. (Note that we follow the convention of \cite{Tu86,Tu01}, the torsion as in \cite{Mi66} is the multiplicative inverse of our torsion.) If $X$ is the empty set, then we write of course $\tau(N,\gamma)$ instead of $\tau(N,X,\gamma)$.

It follows from standard arguments (cf. the above literature) that
the Reidemeister torsion $\tau(N,\gamma)$ is well--defined
up to multiplication by an element of the form $\pm r$ where $r\in \det(\gamma(\pi))$.
Put differently, up to that indeterminacy $\tau(N,\gamma)$ is independent of the choice of underlying CW--structure, the ordering of the cells
and the choice of the lifts of the cells.

Note that $\gamma$ extends to a map $\gamma\colon \Z[\pi]\to M(k\times k,Q)$.
Given an $r\times s$-matrix $A$ over $\Z[\pi]$ we denote by $A_{\gamma}$
 the $rk\times sk$-matrix which is given by applying $\gamma$ to each entry of $A$.
If $B$ is the matrix over $\Z[\pi]$ which represents the boundary map of
$C_i(\ti{N},\ti{X})\to C_{i-1}(\ti{N},\ti{X})$ with respect to the bases given by the lifts,
then $B_{\gamma}$
represents the boundary map of
\[ C_i(\ti{N},\ti{X})\otimes_{\Z[\pi]}Q^k\to C_{i-1}(\ti{N},\ti{X})\otimes_{\Z[\pi]}Q^k\]
 with respect to the aforementioned bases.

\subsection{Twisted Reidemeister torsion of manifold pairs}\label{section:deftorsion2}

Let $(N,\phi)$ be a 3--manifold pair and let $\a\colon \pi_1(N)\to \gl(k,R)$ be a representation over a domain $R$.
We denote by $Q$ the quotient field of $R$.
The representation $\a$ and $\phi\colon \pi_1(N)\to \Z$  then give rise to  a tensor representation
\[ \begin{array}{rcl} \a\otimes \phi\colon  \pi &\to & \mbox{GL}(k,Q(t)) \\
g&\mapsto & \a(g)\cdot t^{\phi(g)}\end{array}. \]
Note that in this case  $\tau(N,\phi\otimes \a)\in Q(t)$ is well-defined up to multiplication by an element in $Q(t)$ of the form
$\pm rt^k$ where $r\in \det(\a(\pi_1(N)))$ and $k\in \Z$.
In particular if $\a$ is a special linear representation, then $\tau(N,\phi\otimes \a)$ is well-defined up to multiplication by an element of the form $\pm t^k, k\in \Z$. Henceforth, when we give an equality for Reidemeister torsion we mean that there exists a representative for which the equality holds. Similarly, when we say that $\tau(N,\phi\otimes \a)$ is monic, then  we mean that there exists a representative which is monic.

The twisted Reidemeister torsion corresponding a to 3--manifold pair and a representation $\a$
was first studied, with somewhat different definitions, by Lin \cite{Li01}, Wada \cite{Wa94} and Kitano \cite{Ki96}.
We refer to the survey paper \cite{FV10} for more information.

\subsection{Turaev's theorem}

We will several times make use of the following theorem, which is easily seen to be a special case of \cite[Theorem~2.2]{Tu01}.

\begin{theorem}\label{thm:turaev}
Let $Q$ be a field and let
\[C_* \,\,\,= \,\,\, 0\to Q^{n_3}\xrightarrow{B_{3}} Q^{n_2}\xrightarrow{B_2} Q^{n_1} \xrightarrow{B_{1}}Q^{n_0}\to  0 \]
be a complex.
We pick a subset of rows from $B_3$ and a subset of columns from $B_1$
and we delete the corresponding columns and rows from $B_2$
in such a way that we obtain square matrices $A_3,A_2$ and $A_1$.
If $\det(A_3)\ne 0$ and $\det(A_1)\ne 0$, then
\[ \tau(C_*)=\det(A_3)^{-1}\cdot \det(A_2)\cdot \det(A_1)^{-1}.\]
\end{theorem}

\section{Proof of Theorem \ref{mainthm1}}

Let $(N,\phi)$ be a fibered 3--manifold pair with $N\ne S^1\times D^2$ and $N\ne S^1\times S^2$.
Let $\a\colon \pi_1(N)\to \gl(k,R)$ be a representation over a domain $R$.
We denote by $\S$ the fiber of the fibration and we denote
by $f\colon \S\to \S$ the monodromy. Note that our restriction on $N$ implies that $\S\ne D^2$ and $\S\ne S^2$.
Since $\phi$ is primitive  it follows that $\S$ is furthermore connected.

We henceforth identify $N$ with $(\S\times [0,2])/(x,0)\sim (f(x),2)$ and we identify $\S$ with $\S\times 0$.
We pick once and for all a base point $P$ for $N$ in $\S\times (1,2)$. We furthermore denote by $\ti{N}$ the universal cover of $N$,
which we identify with the set of homotopy classes of paths starting at the base point. We write $\pi=\pi_1(N,P)$ and $\G:=\pi_1(\S\times [1,2])$.
We also pick a curve $\mu$ based at $P$ which intersects $\S$ precisely once and such that the intersection is positive. Note that $\phi(\mu)=1$.

We now endow $\S$ with a CW--structure with exactly one $0$-cell $d_0$ and exactly one $2$-cell $d_{2}$. We denote by $d_{11},\dots,d_{1n}$ the $1$--cells of $\S$.
We can then endow $N=(\S\times [0,2])/\sim$ with a CW--structure by extending the product CW-structure
on $\S\times [0,1]$ to a CW-structure on $N$. Note that we can extend the CW-structure such that there
are no $0$-cells in $\S\times (1,2)$ and such that there is precisely one $3$-cell in $\S\times (1,2)$.
(But note that in general one can not arrange the CW--structure on $\S\times [1,2]$ to be again a product structure.)
Also note that we can arrange that there exists a  $1$-cell $f$ of $\S\times (1,2)$ such that $\partial f=d_0\cup -\ol{d}_0$.
Summarizing, we can endow $N=(\S\times [0,2])/(x,0)\sim (f(x),2)$ with a CW--structure where the cells are given as follows:
\bn
\item $d_{0}:=d_{0}\times \{0\}$ and $\ol{d}_{0}:=d_0\times \{1\}$,
\item $d_{1j}:=d_{1j}\times \{0\}$ and $\ol{d}_{1j}:=d_{1j}\times \{1\}$ for $j=1,\dots,n$,
\item $d_{2}:=d_{2}\times \{0\}$ and $\ol{d}_{2}:=d_2\times \{1\}$,
\item $e_1:=d_0\times (0,1)$,
\item $e_{2j}:=d_{1j}\times (0,1)$ for $j=1,\dots,n$,
\item $e_3:=d_2\times (0,1)$,
\item one $1$-cell $f_1$ with $\partial \ol{f_{1}}=d_0\cup -\ol{d}_0$
\item one $3$-cell $f_{3}$ in $\S\times (1,2)$,
\en
together with a collection $F_1$ of 1-cells in $M:=\S\times (1,2)$ and a collection $F_2$ of $2$-cells in $M$.
\begin{figure}
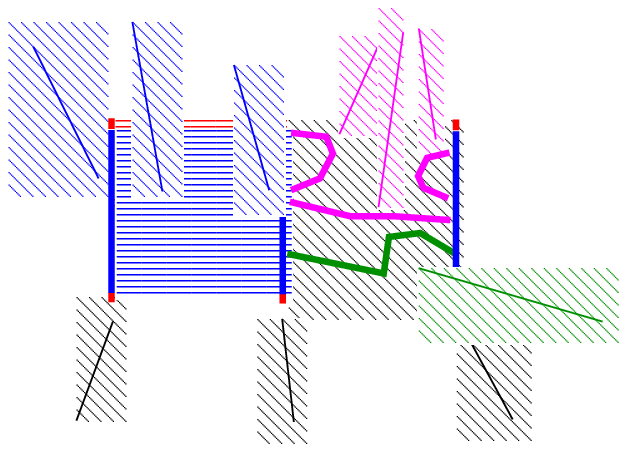
\end{figure}
For each cell we now pick a base point.
Furthermore, for each cell in $\S\times [1,2]$ we  pick a path in $\S\times [1,2]$ from the base point $P$ to the chosen base points in the cells.
We also pick paths in $\S\times (0,2]$ from the base point $P$ to each cell in $\S\times (0,1)$ which intersect $\S=\S\times 0$ precisely once.   Note that these paths define lifts of the cells to $\ti{N}$.
By a slight abuse of notation we denote the lifts by the same symbols.
Note that we can and will pick the orientation of our cells and the basings of our cells such that
\[ \ba{rcl}   \partial e_3&=&d_2-\mu \ol{d}_2\\
            \partial e_{2j}&=&d_{1j}-\mu \ol{d}_{1j}, \mbox{ for } j=1,\dots,n,\\
            \partial e_1&=&d_0-\mu \ol{d}_0,\\
            \partial f_3&=&d_2-z\ol{d}_2+\mbox{linear combination over $\Z[\G]$ of cells in $F_2$},\\
            \partial f_1&=&d_0\cup -x\ol{d}_0
           \ea \]
for some $z,x\in \G$.
We now write
\[  D_1:=\{{d}_{11},\dots,{d}_{1n}\}, \ol{D}_1:=\{{\ol{d}}_{11},\dots,{\ol{d}}_{1n}\}\mbox{ and }E_{2}:=\{{e}_{21},\dots,{e}_{2n}\}.\]
For $i=3,2,1,0$ we now equip
the free $\Z[\pi]$-modules $C_i(\ti{N})$ with the bases
\[ \{e_3,f_3\}, \{E_2,d_2,\ol{d}_2,F_2\},\{e_1,D_1,\ol{D}_1,f_1,F_1\}\mbox{ and } \{d_0,\ol{d}_0\}.\]
Note that with respect to these bases the chain complex is now of the following form
\[ 0\to C_3(\ti{N})\xrightarrow{\bp *&0 \\ 1&1 \\ -\mu&-z\\ 0&*\ep} C_2(\ti{N})
\xrightarrow{\bp *&0&0&0 \\ \id_n&*&0&A \\ -\mu\id_n&0&*&\ol{A} \\ 0&0&0&*\\ 0&0&0&B\ep}
C_1(\ti{N})\xrightarrow{\bp 1 &*&0&1&* \\ -\mu &0&*&-x&* \ep}C_0(\ti{N})\to 0.\]
Here we view the boundary matrices as block matrices corresponding in an obvious fashion to the blocks of basis vectors.
Note that $A,\ol{A}$ and $B$ are matrices with entries in $\Z[\G]$.
\medskip

We now tensor this chain complex with the $\Z[\pi]$-module $Q(t)^k$.
As discussed, the boundary matrices are then given by applying $\a\otimes \phi$ to the above boundary matrices.
We pick the rows of $\partial_3$ corresponding to $d_2\otimes v_1,\dots,d_2\otimes v_k$ and  $\ol{d}_2\otimes v_1,\dots,\ol{d}_2\otimes v_k$
and we pick the columns of $\partial_1$ corresponding to $e_1\otimes v_1,\dots,e_1\otimes v_k$ and $f_1\otimes v_1,\dots,f_1\otimes v_k$.
It now follows from
Theorem \ref{thm:turaev} that
\[\tau(N,\phi\otimes \a)=\det\bp 1&1 \\ -\mu&-z\ep_{\phi\otimes\a}^{-1}\det\bp  \id_n &A\\ -\mu\id_n& \ol{A} \\ 0&B \ep_{\phi\otimes\a} \det\bp 1&1 \\ -\mu&-x\ep_{\phi\otimes\a}^{-1}.\]
Note that $(\phi\otimes\a)(\mu)=t\a(\mu)$ and that $\phi$ vanishes on $\G$. We thus obtain the following equality:
\be \label{equ:tau1} \ba{cl}&\tau(N,\phi\otimes \a)\\
=&\hspace{-0.2cm}\det \bp \id_k&\id_k \\ -t\a(\mu)&-\a(z)\ep^{-1} \hspace{-0.2cm}
\det\bp  \id_{nk} &A_\a\\ (-t\mu \id_n)_{\a}& \ol{A}_{\a} \\ 0&B_\a \ep
\det \bp \id_k&\id_k \\ -t\a(\mu_k)&-\a(x)\ep^{-1}\\
=&\hspace{-0.2cm}\det(t\a(\mu)-\a(z))^{-1}\hspace{-0.1cm}
\left((-t)^{nk}\det(\a(\mu))^n\det\hspace{-0.1cm}\bp A\\ B\ep_{\hspace{-0.1cm}\a}\hspace{-0.2cm}+\dots+\det\hspace{-0.1cm}\bp \ol{A}\\ B\ep_{\hspace{-0.1cm}\a} \right)
\det(t\a(\mu)-\a(x))^{-1}.\ea \ee
(Here and throughout the rest of the paper all calculations will be performed up to sign.)
We will now prove the following claim.

\begin{claim}
There exist $\ol{g},g\in \pi$ such that
\[ \det\bp \ol{A}\\ B\ep_\a =\det(\a(\ol{g})) \mbox{ and }  \det\bp A\\ B\ep_\a =\det(\a(g)).\]
\end{claim}

We  identify $\S$ with $\S\times 0=\S\times 2$ and equip $\ol{M}=\S\times [1,2]$ with the base point $P$.
We denote by $p\colon \widehat{\ol{M}}\to \ol{M}$ the universal covering of $\ol{M}$ and we write $\widehat{\S}:=p^{-1}(\S)$.
 Note that the cells $d_{ij},\ol{d}_{ij}$ and $f_{ij}$ in $\ti{N}$ are in fact naturally cells in $\widehat{\ol{M}}$.
For $i=3,2,1,0$ we now equip $C_*(\widehat{\ol{M}},\widehat{\S})$ with the bases
\[ \{f_3\}, \{\ol{d}_2,F_2\}, \{\ol{D}_1,f_1,F_1\} \mbox{ and } \{\ol{d}_0\}.\]
It follows from the above that the chain complex $C_*(\widehat{\ol{M}},\widehat{\S})$
with the above bases is of the form
\[ 0\to C_3(\widehat{\ol{M}},\widehat{\S})\xrightarrow{\bp -z \\ *\ep}C_2(\widehat{\ol{M}},\widehat{\S})
\xrightarrow{\bp *&\ol{A}\\ *&*\\ *&B\ep}C_1(\widehat{\ol{M}},\widehat{\S})
\xrightarrow{\bp *&-x&*\ep}C_0(\widehat{\ol{M}},\widehat{\S})
\to 0.\]
We again  apply $\a$ to the  boundary matrices.
We then pick the  rows of $\partial_3$ corresponding to   $\ol{d}_2\otimes v_1,\dots,\ol{d}_2\otimes v_k$
and we pick the columns of $\partial_1$ corresponding to  $f_1\otimes v_1,\dots,f_1\otimes v_k$.
It now follows from
Theorem \ref{thm:turaev} that
\[\tau(\ol{M},\S,\a)=\det(\a(z))^{-1}\cdot \det\bp \ol{A}\\ B\ep_\a\cdot \det(\a(x))^{-1}.\]
On the other hand the inclusion map $S\to M$ is a homotopy equivalence.
Since the Whitehead group of a surface group is trivial (see e.g. \cite[p~.250]{Wal78}) this implies by \cite{Mi66} that
the relative torsion is trivial for any coefficient system, i.e. $\tau(M,S,\phi,\a)=1$. We now see that $\ol{g}=xz$ has the desired property.
This concludes the proof of the first statement of the claim. The claim regarding the second matrix is proved exactly the same way.

We now return to the proof of the theorem. Note that the first and the third term in (\ref{equ:tau1})
are monic. The claim now implies that the middle term is also monic.
Together this implies that $\tau(N,\phi\otimes \a)$ is monic.
Furthermore, it follows from (\ref{equ:tau1}) and the above claim that
\[ \deg(N,\phi\otimes\a)=-k+nk-k=k(n-2)=-k\chi_-(\S)=\|\phi\|_T.\]
Here the last equality follows from the well-known fact that a fiber is Thurston norm minimizing.
(In fact this is also an immediate consequence of Theorem \ref{mainthm2}).
This now concludes the proof of Theorem \ref{mainthm1}.

\section{Proof of Theorem \ref{mainthm2}}

Let $(N,\phi)$ be a 3--manifold pair
and let $\a\colon \pi_1(N)\to \gl(k,R)$ be a representation over a domain $R$.
It follows easily from  \cite[Section~1]{Tu02} that we can find a  surface $\S\subset N$ with components $\S_1,\dots,\S_l$ and $r_1,\dots,r_l\in \N$ with the following properties:
\bn
\item $r_1[\S_1]+\dots+r_l[\S_l]$ is dual to $\phi$,
\item $\sum_{i=1}^l -r_i\chi(\S_i)=\|\phi\|_T$,
\item  $N\sm \S$ is connected.
\en
The proof of Theorem \ref{mainthm2} proceeds in a similar fashion to the proof of Theorem \ref{mainthm1}
by picking a suitable CW-structure. Since the surface is now disconnected the notation becomes necessarily more heavy.

For $i=1,\dots,l$  we pick disjoint oriented tubular neighborhoods $\S_i\times [-1,2]$ and we identify $\S_i$ with $\S_i\times \{0\}$.
We write $M:=N\sm \cup_{i=1}^l \S_i\times [0,1]$.
We pick once and for all a base point $P$ in $M$ and we denote by $\ti{N}$ the universal cover of $N$.
We write $\pi=\pi_1(N,P)$ and $\G:=\pi_1(M,P)$.
For $i=1,\dots,l$ we also pick a curve $\mu_i$ based at $P$ which intersects $\S_i$ precisely once in a positive direction
and does not intersect any other component of $\S$. Note that $\phi(\mu_i)=r_i$.

We now build a CW--structure on $N$ as follows.
For each $i=1,\dots,l$ we first endow $\S_i$ with a CW--structure with exactly one $0$-cell $d_0^i$, exactly one $2$-cell $d_{2}^i$
 and $1$-cells $d_{11}^i,\dots,d_{1n_i}^i$.
For $i=1,\dots,l$ we then equip $\S_i\times [-1,0]$, $\S_i\times [0,1]$ and $\S_i\times [1,2]$
with product CW-structures. Since $M$ is connected we can pick $l$ disjoint curves
which connect a point  in $d_2^i\times -1$ with a point in $d_2^i\times 2$.
We use these curves to tube the $3$-cells $d_2^i\times (-1,0)$ and $d_2^i\times (1,2)$.
We denote the resulting $3$-cells by $f_3^1,\dots,f_3^l$. We then extend the CW-structure to a CW-structure on all of $N$.
Since $M$ is connected  we can arrange that there are no $0$-cells in $M$.
Furthermore, by `swallowing' other $3$-cells we can in fact arrange that $f_3^1,\dots,f_3^l$ are the only $3$-cells in $M$.
Finally we  can arrange that
 for $i=1,\dots,l$ there exists a $1$-cell $f_1^i$  such that $\partial f_1^i=
d_0^i\cup -\ol{d}^i_0$.
Summarizing, we can endow $N$ with a CW--structure where for $i=1,\dots,l$ we have the following cells:
\bn
\item $d^i_{0}:=d^i_{0}\times \{0\}$ and $\ol{d}^i_{0}:=d^i_0\times \{1\}$,
\item $d^i_{1j}:=d^i_{1j}\times \{0\}$ and $\ol{d}^i_{1j}:=d^i_{1j}\times \{1\}$ for $j=1,\dots,n_i$,
\item $d^i_{2}:=d^i_{2}\times \{0\}$ and $\ol{d}^i_{2}:=d^i_2\times \{1\}$,
\item $e^i_1:=d^i_0\times (0,1)$,
\item $e^i_{2j}:=d^i_{1j}\times (0,1)$ for $j=1,\dots,n_i$,
\item $e^i_3:=d^i_2\times (0,1)$,
\item one $1$-cell $f_1^i$ in $M$ with $\partial {f_{1}^i}=d^i_0\cup -\ol{d}^i_0$,
\item one $3$-cell $f_3^i$ in $M$ with
\[  \partial f_3=\sum_{i=1}^kd_2^i-\ol{d}_2^i+\mbox{linear combination of cells in $M$.}\]
\en
and there is a collection $F_1$ of $1$-cells in $M$ and a collection $F_2$ of $2$-cells in $M$.
For each cell we now pick a base point and for each cell in $\ol{M}$ we pick a path in $\ol{M}$ from the base point $P$ to the chosen base points.
Furthermore for each cell in $\S_i\times (0,1)$ we pick a path in $M\cup \S_i\times (0,2]$ from the cell to  the base point $P$.
These paths define lifts of the cells to $\ti{N}$ and by a slight abuse of notation we denote the lifts again by the same symbols.
Note that we can and will pick the orientation of our cells and the basings of our cells such that for $i=1,\dots,l$ we have
\[ \ba{rcl} \partial e^i_3&=&d^i_2-\mu \ol{d}^i_2\\
            \partial e^i_{2j}&=&d^i_{1j}-\mu \ol{d}^i_{1j}\mbox{ for } j=1,\dots,n_i,\\
            \partial e^i_1&=&d^i_0-\mu \ol{d}^i_0\ea \]
and such that
\[ \ba{rcl}  \partial f_3^i&=&\sum_{i=1}^kd_2^i-z_i\ol{d}_2^i+\mbox{linear combination over $\Z[\G]$ of cells in $F_2$},\\
            \partial f_1^i&=&d_0^i\cup -x_i\ol{d}_0^i,\ea \]
where $x_1,\dots x_l$ and $z_1,\dots,z_l$ lie in $\G$.
For $i=1,\dots,l$ we write
\[ D^i_1:=\{d^i_{11},\dots,d^i_{1n}\}, \ol{D}^i_1:=\{\ol{d}^i_{11},\dots,\ol{d}^i_{1n}\}\mbox{ and }E^i_{2}:=\{e^i_{21},\dots,e^i_{2n}\}.\]

In the remaining discussion we now only consider the case $l=2$ to simplify the notation.
It should be obvious to the reader that the general case can be treated exactly the same way.

 Note that the chain groups  $C_i(\ti{N})$  are free $\Z[\pi]$-modules. For $i=3,2,1,0$ we now equip them with the bases
\[ \{e_3^1,f_3^1,e_3^2,f_3^2\}, \{E_2^1,E_2^2,d_2^1,\ol{d}_2^1,d_2^2,\ol{d}_2^2,F_2\},\{e_1^1,e_1^2,D_1^1,\ol{D}_1^1,D_1^2,\ol{D}_1^2,f_1^1,f_1^2,F_1\}\mbox{ and } \{d_0^1,\ol{d}_0^1,d_0^2,\ol{d}_0^2\}.\]
With respect to these bases the boundary maps are then  given by the following matrices:
\[ \ba{rcl} B_3&=&\bp *&0&0&0 \\0&0&*&0 \\ 1&1&0&0 \\ -\mu_1&-z_1&0&0 \\0&0&1&1 \\ 0&0&-\mu_1&-z_2\\ 0&*&0&* \ep\\
B_2&=&\bp *&0&0&0&0&0&0 \\  0&*&0&0&0&0&0 \\
\id_{n_1}& 0 &*&0&0&0&* \\
-\mu_1\id_{n_1}&0&0&*&0&0&*\\
0&\id_{n_2}& 0&0&*&0&* \\
0&-\mu_2\id_{n_2}&0&0&0&*&*\\
0&0&0&0&0&0&* \\
0&0&0&0&0&0&* \\
0&0&0&0&0&0&*\ep \\
B_3&=&
\bp 1& 0&*&0&0&0& 1&0&* \\
-\mu_1&0&0&*&0&0&-x_1&0&* \\
0&1&0&0&*&0&0&1&* \\
0&-\mu_2 &0&0&0&*&0&-x_2&* \ep \ea \]
Here we view the boundary matrices as block matrices corresponding in an obvious fashion to the blocks of basis vectors.
Note that all matrices marked by $*$ are matrices with entries in $\Z[\G]$.
\medskip

We now tensor this chain complex with the $\Z[\pi]$-module $Q(t)^k$.
The boundary matrices are then given by applying $\a\otimes \phi$ to the above boundary matrices.
We pick the rows of $\partial_3$ corresponding to $d_2^i\otimes v_j$ and $\ol{d}_2^i\otimes v_j$ for $i=1,2$, $j=1,\dots,k$
and we pick the columns of $\partial_1$ corresponding to $e_1^i\otimes v_j$ and $f_1^i\otimes v_j$ for $i=1,2$, $j=1,\dots,k$.
It now follows from
Theorem \ref{thm:turaev} that $\tau(N,\phi\otimes \a)$ equals
\[
\det\bp 1&1&0&0 \\ -\mu_1&-z_1&0&0 \\ 0&0&1&1 \\ 0&0&-\mu_2&-z_2\ep_{\hspace{-0.1cm}\phi\otimes\a}^{-1}
\hspace{-0.2cm}
\det\bp \id_{n_1}& 0 &* \\ -\mu_1\id_{n_1}&0&*\\ 0&\id_{n_2}& * \\ 0&-\mu_2\id_{n_2}&*\ep_{\hspace{-0.1cm}\phi\otimes\a}
\hspace{-0.2cm}
\det\bp 1&0&1&0 \\ -\mu_1&0&-x_1&0 \\ 0&1&0&1 \\ 0&-\mu_2&0&-x_2\ep_{\hspace{-0.1cm}\phi\otimes\a}^{-1}.\]
Note that $(\phi\otimes \a)(\mu_i)=t^{r_i}\a(\mu_i)$. It follows that
\[ \ba{rcl} \det\bp 1&1&0&0 \\ -\mu_1&-z_1&0&0 \\ 0&0&1&1 \\ 0&0&\mu_2&-z_2\ep_{\phi\otimes\a}
&=&\det(t^{r_1}\a(\mu_1)-\a(z_1))\cdot \det(t^{r_1}\a(\mu_2)-\a(z_2))\\
&=&\left(t^{kr_1}\det(\a(\mu_1))+\dots\pm\a(z_1)\right)\cdot \left(t^{kr_2}\det(\a(\mu_2))+\dots\pm\a(z_2)\right)
\ea
\]
is a polynomial of degree $kr_1+kr_2$. The same argument shows that the degree of the determinant of the third matrix in the above calculation of $\tau(N,\phi\otimes \a)$ equals $kr_1+kr_2$.
Also  note that if we apply the representation $\phi\otimes\a$ to a matrix over $\Z[\G]$
we obtain a matrix with entries in $Q$. Combining these observations we see that there exists a matrix $A$ over $Q$ such that
\[  \deg \tau(N,\phi\otimes \a)=-2kr_1-2kr_2+\deg \det
\left(
\bp 0& 0 &0 \\ t^{r_1}P_1&0&0\\ 0&0& 0 \\ 0&t^{r_2}P_2&0\ep
+A \right)\]
where $P_i=(-\mu_i\id_{n_i})_\a$.
Note that $P_i$ is a $kn_i\times kn_i$-matrix over $Q$.
 It is now elementary to see that
\[ \deg\det
\left(
\bp 0& 0 &0 \\ t^{r_1}P_1&0&0\\ 0&0& 0 \\ 0&t^{r_2}P_2&0\ep
+A \right)\leq kn_1r_1+kn_2r_2.\]
It follows that
\[ \ba{rcl} \deg(\tau(N,\phi\otimes \a))&\leq& -2kr_1-2kr_2+kn_1r_1+kn_2r_2=k((n_1-2)r_1+(n_2-2)r_2)\\
&=&k(r_1\chi(\S_1)+r_2\chi(\S_2))=k\|\phi\|_T.\ea \]
This concludes the proof of Theorem \ref{mainthm2}.

\end{document}